%
%
%
%
%
%
%

\input amstex.tex
\input amsppt.sty
\documentstyle{amsppt}
\advance\vsize2\baselineskip

\def\deg{\mathop{\fam0 deg}}
\define\Z{\Bbb Z}
\define\R{\Bbb R} 
\define\rest#1{\left|_{#1}\right.}
\rightheadtext{}
\leftheadtext{}
\NoBlackBoxes

\topmatter
\title Classification of framed links in 3-manifolds \endtitle
\author
Matija Cencelj, Du\v san Repov\v s and Mikhail B. Skopenkov
\endauthor
\address
Institute for Mathematics, Physics and Mechanics, and
Faculty of Education, University of Ljubljana, P. O. Box 2964, 1001 Ljubljana,
Slovenia.  \endaddress
\email
matija.cencelj\@guest.arnes.si \endemail
\address Institute for Mathematics, Physics and Mechanics, and
Faculty of Education, University of Ljubljana, P. O. Box 2964, 1001 Ljubljana,
Slovenia.  \endaddress
\email
dusan.repovs\@guest.arnes.si \endemail
\address
Department of Differential Geometry, Faculty of Mechanics and
Mathematics, Moscow State University, Moscow, Russia,
119992.\endaddress \email skopenkov\@rambler.ru \endemail
\abstract We present a short complete proof of the following Pontryagin
theorem, whose original proof was complicated and has never been
published in details: {\sl Let $M$ be a connected
oriented closed smooth 3-manifold, $L_1(M)$ be the set of
framed links in $M$ up to a framed cobordism, and $\deg:L_1(M)\to
H_1(M;\Z)$ be the map taking a framed link to its homology class.
Then for each $\alpha\in H_1(M;\Z)$ there is a 1-1 correspondence
between the set $\deg\nolimits^{-1}\alpha$ and the group $\Bbb
Z_{2d(\alpha)}$, where $d(\alpha)$ is the divisibility of the
projection of $\alpha$ to the free part of $H_1(M;\Bbb Z)$.}
\endabstract
\endtopmatter

\document

\bigskip
\noindent {\bf  1 Introduction}
\bigskip

Throughout this paper let $M$ be a connected
oriented closed smooth 3-di\-men\-sional manifold.
Denote by $L_1(M)$ the set of 1-di\-men\-sional framed links in $M$
up to a framed cobordism.
The main purpose of this paper is to describe the set $L_1(M)$.

This classification problem appeared in Pontryagin's
investigations connected with the calculation of the homotopy
groups of spheres and, in a more general situation, of cohomotopy
sets. {\it The cohomotopy set} $\pi^2(M)=[M;S^2]$ is the set of
continuous maps $M\to S^2$ up to a homotopy. 
By the
Pontryagin-Thom
construction this set is in 1--1 correspondence with the set
$L_1(M)$.
Notice that the set of all 
nonzero vector fields
on $M$ up to homotopy, as well as
the set of all
oriented plane fields on $M$ up to homotopy, is also in 1--1
correspondence with the set $L_1(M)$, because 
every orientable 3-manifold
is paralellizable.

To state the main result we need the notions of {\it the natural
orientation} on a framed link and {\it the degree} of a framed
link, defined as follows. The link $L$ is {\it naturally oriented}
if for each point $x\in L$ the tangent vector of the orientation
together with the two vectors of the framing gives a positive
basis of $M$. The {\it degree} $\deg L$ of $L$ is the homology
class (with integral coefficients) of naturally oriented $L$. So
we have a map
$$
\deg:L_1(M)\to H_1(M;\Z).
$$
The classical Hopf-Whitney theorem (1932-35) asserts that this map
is always surjective.

\proclaim{Theorem 1 (L.S.~Pontryagin)} Let $M^3$ be a connected
oriented closed smooth 3-manifold. Then for each $\alpha\in
H_1(M^3;\Z)$ there is a 1-1 correspondence between the sets
$\deg\nolimits^{-1}\alpha$
and
$\Bbb Z_{2d(\alpha)},$
where $d(\alpha)$ is the divisibility of the projection of
$\alpha$ to the free part of $H_1(M^3;\Bbb Z)$.
\endproclaim

\proclaim{Example} The set of all maps $f:S^1\times S^1\times S^1\to
S^2$ is up to homotopy in a
bijective correspondence with the set of all
$4$-tuples
$(p,q,r,t)$, where $p,q,r\in\Bbb Z$ are the degrees of the
restrictions of $f$ to the 
2-dimensional subtori, $t\in\Bbb Z$, for
$p=q=r=0$, and $t\in \Bbb Z_{2\gcd(p,q,r)}$, otherwise.
\endproclaim

Recall that the divisibility of zero is zero and the divisibility
of a nonzero element $\alpha\in G$ is $\max\{\,d\in\Bbb
Z\,|\,\exists\beta\in G:\alpha=d\beta\,\}$. We denote $\Bbb
Z_d=\Bbb Z/d\Bbb Z$, in particular, $\Bbb Z_0=\Bbb Z$.

In this paper we give a short and direct proof of
Theorem 1 (which was
stated without proof in [Po]). In fact, Theorem 1 was not even properly
stated
in the paper [Po] itself (the paper was written in English), but only in the abstract (written
in Russian) without any indication of the proof. The statement in
the abstract is not clear, so we have borrowed it from [St].

The statement from [St] asserts that there is a 1-1
correspondence between
$$
\deg\nolimits^{-1}\alpha\qquad\text{ and }\qquad\frac{\Z}{2\alpha\cap H_2(M;\Z)},
$$
which by the Poincar\'e
duality, is equivalent to our statement of  Theorem 1. 
There are reasons to believe that our proof is the same as the
proof which Pontryagin had in mind, but never published it, going
instead straight to the general case when $M$ is an
arbitrary polyhedron.

In the stable codimension $n\ge 4$ there is an analogous theorem
describing the set of framed links in $n$-manifolds [RSS]:

\proclaim{Theorem 2 (L.S. Pontryagin, N. Steenrod, W.T. Wu)}
 Let $M$ be a connected oriented closed smooth
$n$-manifold, $n\ge4$. Then the
degree map $\deg:L_1(M)\to H_1(M;\Z)$
is a bijection, if there exists $\beta\in H_2(M,\Z)$ such that
$\rho_2\beta\cdot w_2(M)=1\pmod2$. If such $\beta$ does not exist,
then $\deg$ is a 2-1 map (that is, each $\alpha\in H_1(M;\Z)$ has
exactly two preimages).
\endproclaim

Here $w_2(M)$ is the Stiefel-Whitney class and $\rho_2:H_1(M;\Bbb
Z)\to H_1(M;\Bbb Z_2)$ is the reduction modulo 2. In this paper we
use an extension of the ideas of [RSS]. (Note that there is a
misprint in the statement of this theorem in [RSS; Theorem 1a]).

The results stated above remain of sufficient interest up to the present -- see
[AK], [BP], [Du], [Go], [GG], and [Ka]. Note that
Proposition~4.1 of [Go] is equivalent to our Theorem~1. 
However, even though [Go] gives
the proof of this proposition, it is not written in details ­- in the notations
of our paper (see \S3 below), the statement 
(1) that the invariant $h$ is well-defined 
and the surjectivity of $h$ are indeed verified, whereas the proof of the {\sl injectivity},
which is not evident, is absent from [Go].

Notice also that the statements of this result
in [BP; Theorem 6.2.7] and [Du; Proposition 1] are erroneous,
because of a different definition of the number $d(\alpha)$.
(In these papers $d(\alpha)$ is defined to be 0 if $\alpha$ is a torsion
element, otherwise it is the divisibility of $\alpha$ in $H_1(M^3;\Bbb
Z)$. This is not equivalent to our definition.) 
An alternative approach to Theorem~1, by  different methods, can be found in [AuKa]. A sketch of alternative
proof can be found in [Ku; Proposition 2.1]. A
more general result appears in [Ka; Theorem 2.4.7] and [GG; Proposition 7].'

The plan of the paper is as follows: in \S2 we first recall a nice geometric definition of the 
normal Euler class
and we then prove another Pontryagin's classification theorem (see Remark
after Lemma~3). In \S3 we finally 
prove Theorem~1.

\bigskip
\noindent
{\bf  2 Preliminaries}
\bigskip

We are going to use the following {\it geometric definition} of the
normal Euler class, which is
equivalent to other definitions. Let $M^4$ be
a closed oriented connected 4-manifold. Let $L^2$ be a connected
oriented manifold immersed in $M^4$. Let $\nu(L)$ be the normal
bundle of $L$. Identify $L$ with the zero section of $\nu(L)$. Fix
a natural orientation of $\nu(L)$. Take a general position section
$L'$ of $\nu(L)$. {\it The  normal Euler class $\bar e(L)=\bar
e(\nu(L))\in\Z$} is the difference between the numbers of positive
and negative intersection points of $L$ and $L'$. Further denote
by $X\cap Y$ the difference between the numbers of positive and
negative intersection points of $X$ and $Y$.

We are also going to use the following {\it geometric definition} of the
relative normal Euler class.
Fix an orientation of $M^3\times [0;1]$.
Let $L_1\subset M\times 1$ and $L_2\subset M\times 0$
be a pair of framed links, let $L\subset M\times[0;1]$
be a (unframed) cobordism between them.
Fix a natural orientation of $L$, i.e. an orientation
that induces natural orientations of $L_1$ and $L_2$.
Fix a natural orientation of $\nu(L)$.
The first vector field of the framings of $L_1$ and $L_2$
can be considered as a section of $\partial \nu(L)$.
Let $L'$ be a general position extension of this section
to a section of $\nu(L)$.
{\it The relative normal Euler class $\bar e(L)\in\Z$}
is the difference between the numbers of positive and negative
intersection points of $L$ and $L'$.
If we reverse the orientation of $M^3\times [0;1]$
(and, consequently, of $L$, because $L$ is naturally oriented),
then the sign of the integer $\bar e(L)$ changes.

It can be shown that the class $\bar e(L)$
is the complete obstruction to extension of the framing of $\partial L$
to a framing of $L$.

\proclaim{Lemma 3} Let $L^2$ and $M^4$ be a pair of connected
oriented manifolds ($M$ may have boundary). Suppose that $L$ is
immersed into $M$. Denote by $[L]\in H_2(M;\Z)$ the class of $L$.
Denote by $\sigma$ the difference between the numbers of positive
and negative self-intersection points of $L$. Then
$$
\bar e(L)=[L]\cdot[L]-2\sigma,
$$
where we identify  the group $H_0(M;\Z)$ with $\Z$.
In particular, if $M=N^3\times I$ for some 3-manifold $N^3$,
then $\bar e(L)=-2\sigma$.
\endproclaim

\example{Remark} In particular, this well-known lemma implies
[RSS; Theorem 1.2b].
\endexample

\demo{Proof of Lemma 3}
Let $\pi$ be the natural projection of a neighbourhood of $L$ in $\nu(L)$
to a small neighbourhood of $L$ in $M$.
Take a general position section $L'$ of $\nu(L)$ close to zero.
The lemma now follows from
$$
\bar e(L)=L\cap L'=\pi L\cap \pi L'-2\sigma=[L]\cdot[L]-2\sigma.\qed
$$
\enddemo

\bigskip
\noindent {\bf  3 Proof of Theorem 1}
\bigskip

In order to construct a bijection
$h:\deg\nolimits^{-1}\alpha\to\Bbb Z_{2d(\alpha)}$, fix a framed
circle $L_1$ such that $\deg L_1=\alpha$ (clearly, such a circle
exists). Take an arbitrary framed link $L_2$ such that $\deg
L_2=\alpha$. Since $L_1$ and $L_2$ are homologous, it follows that
there is a (not framed) cobordism $L$ between them. By definition,
put $h(L_2)=\bar e(L)\mod 2d(\alpha)$. (One can see that this is
the Hopf invariant if $\alpha=0$ and $L_1$ is null framed
cobordant.)

It will follow from (1) and (2) below that $h$ is well-defined:

(1) $h(L_2)$ does not depend on the choice of $L$; and

(2) if $L_2$ and $L_2'$ are framed cobordant then
$h(L_2)=h(L_2')$.

Let us first prove (2). Assume that $L_1\subset M\times 1$,
$L_2\subset M\times 0$, $L_2'\subset M\times (-1)$, $L\subset
M\times [0,1]$. Let $L'\subset M\times [-1,0]$ be a framed
cobordism between $L_2$ and $L_2'$. By the geometric definition of
the relative normal Euler class it follows that $\bar e(L\cup
L')=\bar e(L)+\bar e(L')$. Since the cobordism $L'$ is framed, it
follows that $\bar e(L')=0$. Thus $\bar e(L\cup L')=\bar e(L)$,
and we obtain the required equality $h(L_2)=h(L_2')$.

Let us prove (1).
Take another general position cobordism $L'$ between $L_1$ and $L_2$.
Assume that $L_2\subset M\times 0$,
two copies of $L_1$ are contained in $M\times (\pm1)$ and
$L, L'\subset M\times[0,1]$.
Let $-L'\subset M\times[-1,0]$ be the cobordism symmetric to $L'$
(we consider the symmetry $x\times t\to x\times(-t)$ on $M\times \Bbb R$).
Take a general position framed circle $-L_1'\subset M$
such that $L_1\cup L_1'$ is framed cobordant to zero, i. e. to
an empty submanifold.
Denote by $\Delta$ the corresponding framed cobordism.
Assume that two copies of $L_1'$ are contained in $M\times(\pm1)$,
and $\Delta\subset[1;+\infty)$.
Let $-\Delta\subset (-\infty,-1]$ be the cobordism, symmetric to $\Delta$.
Denote by
$$
K=(-L')\cup L\cup \Delta\cup (L_1'\times [-1,1])\cup(-\Delta).
$$
By the geometric definition of the relative normal Euler class
we obtain
$$\bar e(K)=\bar e(-L')+\bar e(L)+\bar e(\Delta)+\bar e(L_1'\times [-1,1])
+\bar e(-\Delta).
$$
Here $\Delta$, $-\Delta$ and $L_1'\times [-1,1]$ can be framed,
so $\bar e(\Delta)=\bar e(-\Delta)=\bar e(L_1\times [-1,1])=0$.
Since the symmetry $x\times t\to x\times(-t)$ reverses the orientation
of $M\times [-1;1]$, it follows by the geometric definition of the relative
normal Euler class that $\bar e(-L')=-\bar e(L')$.
Thus $\bar e(K)=\bar e(L)-\bar e(L')$.
Now (1) follows from
$$
\bar e(K)=-2\sigma=2(-L'\cup L)\cap (L_1'\times [-1,1])= 2K\cap
(L_1''\times \R)= 2[pK]\cdot \alpha=0\mod 2d(\alpha).
$$
Here $\sigma$ is the difference between the numbers of positive
and negative self-inter\-sec\-tions of $K$, and the first equality
follows from Lemma 3. The second equality follows from the
construction of $K$. Then, $L_1''\subset M$ is a general position
circle close to $L_1'$ and homologic to it. By general position
$L_1'$ and $L_1''$ are disjoint, so $(-L'\cup L)\cap (L_1'\times
[-1,1])=K\cap (L_1''\times [-1,1])$. 

Since $-\Delta$ is obtained
from $\Delta$ by the symmetry $x\times t\to x\times(-t)$, it
follows that $K\cap (L_1''\times [1,+\infty))=-K\cap (L_1''\times
(-\infty,-1])$, and the third equality follows. Denote by
$p:M\times I\to M$ the projection. Then by general position we
obtain the fourth equality, because the homological class of
$L_1''$ is $\alpha$. The last equality follows from the definition
of $d(\alpha)$. So the proof of (1) is completed.

{\bf Injectivity of $h$.} Let $L_2$ and $L_2'$ be a pair of framed
1-sub\-ma\-ni\-folds such that $h(L_2)=h(L_2')$. Let us prove that
$L_2$ and $L_2'$ are framed cobordant. Assume that $L_2\subset
M\times 1$, $L_1\subset M\times 0$ and $L_2'\subset M\times (-1)$.
Let $L\subset M\times[0,1]$ and $-L'\subset M\times[-1,0]$ be the
cobordisms between $L_1$ and $L_2$, $L_1$ and $L_2'$ respectively.
Since $h(L_2)=h(L_2')$, it follows that $\bar e(L)=-\bar
e(-L')\mod 2d(\alpha)$. Then $\bar e(-L' \cup L)=2d(\alpha)y$ for
some $y\in\Z$.

By the Poincar\'e duality there exists an element $\beta\in
H_2(M;\Z)$ such that $\alpha\cap\beta=d(\alpha)$. Let $K\subset
M\times 0$ be a general position connected submanifold realizing
the class $y\beta$. Notice that $\bar e(K)=0$ by Lemma 3. Denote
by $K'$ the connected sum of $(-L'\cup L)$ and $K$ in $M\times
[-1;1]$. By the geometric definition of the relative normal Euler
class it follows that $\bar e(K')=\bar e(-L'\cup L)+\bar
e(K)=2d(\alpha)y$. 

The difference between the numbers of positive
and negative self-intersection points of the manifold $K'$ is
equal to $y\beta\cap\alpha=d(\alpha)y$. Let $K''$ be a new
cobordism between $L_2$ and $L_2'$ obtained from $K'$ by
elimination of the self-intersection points. (Here we use a move
in a neighbourhood of each self-intersection point analogous to
the move taking the pair of the planes $x=0$, $y=0$ and $z=0$,
$t=0$ to the surface
$$
\cases
x(\tau,\varphi)=\tau\cos \varphi,\\
y(\tau,\varphi)=\tau\sin\varphi,\\
z(\tau,\varphi)=(1-\tau)\cos \varphi,\\
t(\tau,\varphi)=(1-\tau)\sin\varphi;
\endcases
$$
in $\R^4$ with coordinats $(x,y,z,t)$.) By the geometric
definition of the normal Euler class it can be proved easily that
removing of each self-intersection point decreases $\bar e(K')$ by
$\pm2$, depending on the sign of the point (since our move is
local, it suffices to prove it for a closed submanifold $K'$, and
this latter case follows from Lemma~3). So $\bar e(K'')=\bar
e(K')-2d(\alpha)y=0$. Thus $K''$ can be framed, hence $L_2$ and
$L_2'$ are framed cobordant.

{\bf Surjectivity of $h$.} Let us construct a sequence $L_1, L_2,
\dots, L_{2d(\alpha)}$ of framed 1-sub\-ma\-ni\-folds such that
for $j=1,\dots,2d(\alpha)$ we have $h(L_j)=j-1$. Fix a
homeomorphism $L_1\cong S^1$. Denote by $f_1(x)$ the basis vector
of the fixed framing of $L_1$ at the point $x\in S^1$. Take a map
$\varphi:S^1\to SO(2)$ realizing the generator
$\pi_1(SO(2))\cong\Z$. For $j=2,\dots,2d(\alpha)$ define the
framing $f_j$ of $L_1$ by the formula
$f_j(x)=\varphi^{j-1}(x)f_1(x)$. Let $L_j$ be the submanifold
$L_1$ with framing $f_j$. Without loss of generality we may assume
that $h(L_2)\ge0$. 

Let us prove that then $h(L_2)=1$. Then it can
be shown analogously that $h(L_j)=j-1$. Take $L=L_1\times I$. It
suffices to construct a general position normal vector field on
$L$ extending the first field of the framing of $L_1$ and $L_2$
with a unique singular point. The normal bundle to $L$ in $M\times
\R$ is trivial. Identify this bundle with $\R\times \R\times L$
and denote by $p_1,p_2:\R\times \R\times L\to\R$ the projections
to the first and the second multiples respectively. Further denote
by $f_2$ the first vector field of the framing $f_2$. Clearly,
$p_1 f_2(x)$, where $x\in L_2$, has exactly two zeros. Join them
by an arc $A\subset L$. 

Analogously, join by an arc $B$ the pair of
zeros of $p_2 f_2(x)$. Clearly, we can choose the arcs $A$ and $B$
intersecting transversally at a single point. Take a general
position normal vector field $F_1$ on $L$ extending the fields
$p_1 f_2$, $p_1 f_1$ and such that $p_2 F_1=0$, $p_1 F_1\left|_A
\right.=0$. Analogously, extend $p_2 f_2$ and $p_2 f_1$ to a
normal vector field $F_2$ such that $p_1 F_2=0$, $p_2
F_2\left|_B\right.=0$. The sum $F_1+F_2$ with a single zero at the
point $A\cap B$ is the required vector field. \qed

\bigskip\noindent
{\bf  Acknowledgements}
\medskip

\noindent
We thank  A.~B.~Skopenkov for comments and suggestions. This paper first appeared as a preprint [CRS1]
in 2003. After we had posted it also at the ArXives  [CRS2]
in 2007, D.~Auckley, E.~Dufraine and U.~Kaiser kindly brought to our attention  [AK], [BP], [Du], [Go], [GG],  and [Ka].
The first two authors were supported in part by the Slovenian
Research Agency program P1-0292-0101-04 and grants BI-RU/06-07/004,013. 
The third author was supported in
part by Russian Foundation for Basic Research grants
No.~02-01-00014-a, 05-01-00993-a, 06-01-72551-NCNIL-a,
07-01-00648-a, INTAS grant No.~06-1000014-6277, President of the
Russian Federation grant for state support of leading science
schools of the Russian Federation NSh-4578.2006.1 and the Russian
Education and Science Agency program RNP 2.1.1.7988.

\enddocument
\end